\documentclass[12pt]{amsart}
 \usepackage{amssymb,latexsym,amsmath}
\begin{document}
\date{October 25, 2006}
%\title[Superposition]{}
\author{Peter Lindqvist}
\address{Department of Mathematics, Norwegian University of Science and
Technology,  NO-7491 Trondheim, Norway}
 \email{peter.lindqvist@math.ntnu.no}
 \urladdr{http://www.math.ntnu.no/\~{}lqvist}
\author{Juan J. Manfredi}
\title{Hardy's Inequality from a logarithmic Caccioppoli Estimate}
\address{Department of Mathematics, University of Pittsburgh,
Pittsburgh, PA 15260,USA}
 \email{manfredi@pitt.edu}
 \urladdr{http://www.pitt.edu/\~{}manfredi}
 \begin{abstract} We give a simple proof of Hardy's inequality based on the logarithmic Cacciopoli inequality for $p$-superharmonic functions.
 \end{abstract}
%\thanks{This paper was written while the first author was visiting the University of
%Pittsburgh. He wishes to acknowledge the hospitality and the stimulating working
%atmosphere at the Department of Mathematics. The second author was partially supported by NSF %award
% DMS-0500983.}
  \dedicatory{Dedicated to Juha Heinonen (1960-2007)}
 \maketitle
\section{Introduction} The so-called Hardy inequality
$$\sum^{\infty}_{n=1} \ \biggl(\frac{a_{1} + a_{2} + \cdots + a_{n}}{n}\biggr)^{p} \leq \biggl(\frac{p}{p-1}\biggr)^{p} \sum^{\infty}_{n=1} a^{p}_{n},$$
where $a_{n} \geq 0$ and $1 < p < \infty$, was given in 1920 by G. Hardy, c.f. [H].  The sharp constant was fixed by E. Landau, cf. [La].  The integral version for $f(x) \geq 0$ is usually written as
$$\int^{\infty}_{0} \biggl(\frac{F(x)}{x}\biggr)^{p} dx \leq \biggl(\frac{p}{p-1}\biggr)^{p} \int^{\infty}_{0} f(x)^{p} dx$$
where $F(x) = \int^{x}_{0} f(t) dt$.  In a finite interval $[a,b]$ the inequality
$$\int^{b}_{a} \ \bigg| \frac{\zeta(x)}{\delta(x)}\bigg|^{p} dx \leq \biggl(\frac{p}{p-1}\biggr)^{p} \int^{b}_{a} |\zeta'(x)|^{p}dx$$
holds, where $\delta(x) = \min(x-a, b-x)$ and $\zeta$ is a smooth function vanishing at the endpoints $a$ and $b$,  cf. [D, Lemma 1.5.1].  As pointed out in [FHT] the Hardy inequality was known to T. Boggio as early as in 1907,  cf. [B].

We are concerned with the following variant of the Hardy inequality, formulated for a bounded domain $\Omega$ in ${\mathbb R}^{n}$ and involving the distance 
$$\delta(x) = \min\limits_{\xi \in \partial \Omega} | x-\xi |$$
to the boundary.  {\it For} $\zeta \in C^{\infty}_{0}(\Omega)$ {\it the Hardy inequality}
\begin{equation}
\int_{\Omega} \bigg|\frac{\zeta(x)}{\delta(x)}\bigg|^{p} dx \leq C_{\Omega} \int_{\Omega}| \nabla \zeta (x)|^{p} dx
\end{equation}
{\it holds with} $C_{\Omega}$ {\it independent of}  $\zeta$, {\it at least for a smooth domain} $\Omega$.  Our objective is to provide a very simple and natural proof, based on a {\it logarithmic Caccioppoli  estimate}.  The case $p = 2$ involves only the fact that the distance function $\delta$ (or some function of $\delta$) is superharmonic.  Needless to say, the idea to use a differential equation is well-known in this connection, (see [B] and [BFT]) but our approach is, as it were, direct and very easy.  Our use of the logarithmic Caccioppoli estimate seems to be new.

Before proceeding to the proof, we mention some other interesting facts, in passing.  First, the inequality is not valid in all domains. See [Hz].  A necessary and sufficient criterion in terms of $p$-capacities of level sets has been given by V. Mazya in [M, 2.3.3, Corollary, p. 11 ].  Second, the dependence of $p$ is rather puzzling.  For example, in the case of the punctured disk $0 < |x| < 1$, the inequality holds except for $p = n$. The Hardy inequality has \lq\lq the open end property\rq\rq: for a fixed domain the values of $p$ for which it is valid form an open set, cf. \cite{KZ}.
Third, the existence of an extremal function $\zeta$ in $W^{1, p}_{0}(\Omega)$ is equivalent to the strict inequality
$$\biggl(\frac{p}{p-1}\biggr)^{p} < C_{\Omega} < \infty$$
for the best constant.  In a convex domain $\Omega$ the best constant is $(p/(p-1))^{p}$ but it is never attained!  For this fascinating phenomenon we refer to [MMP].

\bigskip
\section{THE LOGARITHMIC CACCIOPPOLI  ESTIMATE}
\smallskip

It is well-known that a bounded superharmonic function $v$ in a domain $\Omega$ satisfies the inequality
\begin{equation}
\int_{\Omega} \langle \nabla v, \nabla \varphi\rangle dx \geq 0\end{equation}
for all non-negative $\varphi \in C^{\infty}_{0}(\Omega)$.  We will need only the case when $v$ is a function of the distance $\delta$.  More generally, we say that a function $v$ belonging to $W^{1,p}(\Omega)$ is $p$-superharmonic in $\Omega$, if
\begin{equation}
\int_{\Omega} \langle |\nabla v|^{p-2} \nabla v, \nabla \varphi\rangle dx \geq 0\end{equation}
for all $\varphi \geq 0, \varphi \in C^{\infty}_{0}(\Omega)$.  Usually, it is also required that $v$ be lower semincontinuous, but since we only need a very restricted class of $v$'s this has no bearing here.  (See [L] for an account on $p$-{\it superharmonic}  functions.)  The next lemma is our starting point.
\bigskip

{\bf Lemma.}  ({\it Caccioppoli}) Suppose that $v$ is $p$-superharmonic in $\Omega$ and that $v > 0$.  Then
\begin{equation}
\int_{\Omega} |\zeta\, \nabla \log v|^{p} dx \leq \biggl(\frac{p}{p-1}\biggr)^p \int_{\Omega} |\nabla \zeta |^{p} dx\end{equation}
for all $\zeta \in C^{\infty}_{0} (\Omega)$.
\bigskip

\noindent
{\bf Proof:}  This well-known estimate follows from (3), when the test function 
$\varphi = |\zeta |^{p} v^{1-p}$ is used.  Indeed, one only has to insert
$$\nabla \varphi = p|\zeta |^{p-2}\zeta v^{1-p} \nabla \zeta - (p-1)|\zeta |^{p} v^{-p} \nabla v$$ and use H\"older's inequality.$\Box$

\bigskip

\section{CONVEX DOMAINS}
\smallskip

Let us turn our attention to the distance function
$$\delta |x| = \min\limits_{\xi \in \partial \Omega} \ |x-\xi |,$$
representing the distance from the point $x$ to the boundary of $\Omega$.  Since
$$|\delta(x) - \delta(y)| \leq |x-y|, \ \ \ x, y\in \Omega,$$
it follows from Rademacher's theorem that $\delta$ is differentiable almost everywhere.  In fact, the eiconal  equation
$$|\nabla \delta (x)|^{2} = 1$$
holds at a.e. $x$ in $\Omega$; see, for example,  Theorem 3.3 in \cite{EH} for this fact.

Notice that if $\Omega$ is a convex domain, then the distance function is superharmonic, because it is the  pointwise minimum of affine functions (hyperplanes).  It is also $p$-superharmonic for any $p$ in the range $1< p < \infty$ (and even for $p = \infty$).  Thus we may use the function $v = \delta$ in the Caccioppoli estimate.  Then
$$|\zeta \nabla \log v | = \frac{|\zeta |}{\delta}$$
and the Hardy inequality below follows.
\bigskip

{\bf Proposition.}  ({\it Hardy}) Let $\Omega$ be a convex domain.  Then
\begin{equation}
\int \bigg|\frac{\zeta(x)}{\delta(x)} \bigg|^{p} dx \leq \biggl(\frac{p}{p-1}\biggr)^p \int_{\Omega} |\nabla \zeta(x)|^{p} dx\end{equation}
holds for all $\zeta \in C^{\infty}_{0}(\Omega).$
\bigskip

{\bf Remark.}  (1) In this case there is no extremal function $\zeta$ in $W^{1,p}_{0}(\Omega)$, although the constant is sharp.  The improved inequality
$$\epsilon \int_{\Omega} |\zeta |^{p} dx + \int_{\Omega} \bigg| \frac{\zeta}{\delta}\bigg|^{p} dx \leq \biggl(\frac{p}{p-1}\biggr)^{p} \int_{\Omega} |\nabla \zeta |^{p} dx$$
is known to be valid for convex bounded domains, c.f.  [MMP].

(2) Since $|\nabla \delta |^{p-2}\nabla \delta = \nabla \delta$ in (3), we can obtain (5) with some effort without referring to $p$-superharmonic functions, using only ordinary superharmonic functions.

\bigskip

\section{THE CASE $p > n$}
\smallskip

The case $p > n$ is as easy as the previous proof.  The functions $|x-\xi |^{(p-n)/(p-1)}$ are $p$-superharmonic when $x\neq \xi$.  We have
$$\delta(x)^{\alpha} = \min_{\xi \in \partial \Omega} | x-\xi |^{\alpha}, \ \ \alpha = \frac{p-n}{p-1} > 0.$$
Hence $\delta^{\alpha}$ is $p$-superharmonic.  Inserting $v= \delta^{\alpha}$ and using
$$|\nabla \log \delta^{\alpha}|^{p} = \alpha^{p}|\nabla \log \delta |^{p}$$
in the Caccioppoli estimate, we again obtain the Hardy inequality, but with the constant
$$\alpha^{-p}\biggl(\frac{p}{p-1}\biggr)^{p} = \biggl(\frac{p}{p-n}\biggr)^{p}.$$
Thus we have proved the following:
\bigskip

{\bf Proposition} ({\it Hardy}) Suppose that $p > n$.  Then the Hardy inequality holds for an arbitrary domain $\Omega$ in ${\mathbb R}^{n},$ with the constant
$$C_{\Omega} =  \biggl(\frac{p}{p-n}\biggr)^{p}.$$

\bigskip

{\bf Remark.}  The constant is sharp in annular domains.

\bigskip

\section{THE CASE $1 < p \leq n$}
\smallskip

For these values of $p$, the Hardy inequality is not valid in arbitrary domains.  Under suitable assumptions the deduction from the logarithmic Caccioppoli estimate works, but at the expense of a longer calculation.

Just to illustrate the method, we assume that $\Omega$ satisfies {\it the exterior sphere condition} with a fixed ball size, say of radius $R$.  We construct the function
$$ v(x) = \min_{y} \biggl(\frac{1}{R^{\beta}} \ - \ \frac{1}{|x-y|^{\beta}}\biggr), \ \ \beta = \frac{n-p}{p-1} > 0$$
where the minimum is taken over all points $y$ satisfying
$$\mbox{dist} \ (y, \overline{\Omega})\geq R.$$
As the pointwise minimum of $p$-harmonic functions $v$ is $p$-superharmonic.

We can write
$$v(x) = \frac{1}{R^{\beta}} \ - \ \frac{1}{(R+\delta(x))^{\beta}}$$
because
$$\min_{y}   |x-y| = R + \delta(x)$$
according to a simple argument involving the exterior sphere and the triangle inequality.

A substitution of the above $v$ into the logarithmic Caccioppoli inequality yields the Hardy inequality.  To this end, we calculate 
$$\nabla v = \ \frac{\beta\nabla \delta}{(R + \delta)^{\beta+1}}.$$
Recalling that $|\nabla \delta | = 1$ a.e., we obtain 
$$\frac{|\nabla v|}{v} \ = \ \beta \ \frac{R^{\beta}}{(R+\delta)[(R + \delta)^{\beta}- R^{\beta}]}$$
The elementary inequality
\[(R+\delta)^{\beta} - R^{\beta} \leq \left\{
\begin{array}{l}
\beta\delta (R + \delta)^{\beta-1}, \ \beta \geq 1\\
\beta\delta R^{\beta-1}, \ 0 < \beta < 1\end{array} \right.\]
yields
\[\frac{|\nabla v|}{v} \ \geq \ \left\{
\begin{array}{l}
\displaystyle{\frac{1}{\delta}} \ \biggl(\frac{R}{R+ \delta}\biggr)^{\beta}, \ \beta \geq 1\\
\displaystyle{\frac{1}{\delta}} \ \frac{R}{R+\delta}, \ 0< \beta < 1\end{array}\right.\]
Therefore,
$$|\nabla\log(v)|^{p} \geq \frac{C}{\delta^p}$$
where $C$ is a power of $R/(R+ \ \mbox{diam} \ \Omega)$.  Again the Hardy inequality follows, since
$$\int\limits_{\Omega} \zeta^{p} | \nabla\log v|^{p} dx \geq C \int\limits_{\Omega} \bigg|\frac{\zeta}{\delta}\bigg|^{p} \ dx.$$
This concludes our proof.

%\section{Heisenberg group}

\end{document}